\documentclass[11pt]{article}

\usepackage{graphicx}
\usepackage{amsmath}
\usepackage{amsthm}
\usepackage{amssymb}

\theoremstyle{plain}
\newtheorem{theorem}{Theorem}

\newtheorem{lemma}[theorem]{Lemma}
\newtheorem{corollary}[theorem]{Corollary}

\theoremstyle{remark}

\newtheorem{remark}[theorem]{Remark}

\newcommand\Mod{{\rm Mod}}
\newcommand\Sym{{\rm Sym}}

\title{Generating Mapping Class Groups by Involutions}
\author{Martin Kassabov\thanks{The author is supported by several NSERC grants}}

\begin{document}
\maketitle
\begin{abstract}
Let $\Sigma_{g,b}$ denote a closed oriented surface genus $g$ with
$b$ punctures and let $\Mod_{g,b}$ denote its mapping class group.
In~\cite{Luo} Luo proved that if the genus is at least $3$, the
group $\Mod_{g,b}$ is generated by involutions. He also asked if
there exists a universal upper bound, independent of genus and the
number of punctures, for the number of torsion
elements/involutions  needed to generate $\Mod_{g,b}$. Brendle and
Farb in~\cite{BF} gave a partial answer in the case of closed
surfaces and surfaces with one puncture, by describing a
generating set consisting of $7$ involutions. Our main result
generalizes the above result to the case of multiple punctures. We also
show that the mapping class group can be generated by smaller
number of involutions. More precisely, we prove that the mapping
class group can be generated by $4$ involutions if the genus $g$
is large enough. There is not a lot room to improve this bound
because to generate this group we need at lest $3$ involutions. In
the case of small genus (but at least $3$) to generate the whole
mapping class group we need a few more involutions.
\end{abstract}

\begin{section}{Introduction}

Let $\Sigma_{g,b}$ denote a closed, oriented surface of genus $g$
with $b$ punctures. By $\Mod_{g,b}$ we will denote its mapping
class group --- the group orientation-preserving homeomorphisms
preserving the set of punctures (i.e., we allow
homeomorphisms of the surface which permute the punctures) modulo
homotopy. Let us also denote by $\Mod^0_{g,b}$ the subgroup of
$\Mod_{g,b}$, which fixes the punctures pointwise. It is clear
that we have the exact sequence:
$$
1 \to \Mod^0_{g,b} \to \Mod_{g,b} \to \Sym_b,
$$
where the last projection is given by the restriction of a homeomorphism
to its action on the puncture points.

The study of the mapping class group of a closed surface begun in
the 1930-es. The first generating set of the mapping class group
$\Mod_{g,0}$ for $g\geq 3$ was constructed by Dehn (see
\cite{De}). This set consists of $2g(g-1)$ Dehn twists. Thirty
years latter, Lickorish (see \cite{Li}) constructed a generating
set for $\Mod_{g,0}$ consisting of $3g - 1$ twists for any $g \geq
1$. This result was improved by Humphries (see \cite{Hu}), who
showed that a certain subset of Lickorish's set consisting of $2g
+ 1$ twists suffices to generate $\Mod_{g,0}$, and that this is in
fact the minimal number of twist generators. This result was
generalized in~\cite{Jo} by Johnson who proved that the same set
of Dehn twists also generates $\Mod_{g,1}$. In the case of
multiple punctures the mapping class group can
be generated by $2g+b$ twists for $b\geq 1$ (see~\cite{Ge}).

It is possible to obtain smaller generating sets of $\Mod_{g,b}$
by using elements other than twists. N.~Lu (see~\cite{Lu})
constructed a generated set of $\Mod_{g,0}$ consisting  of 3 elements
two of which are torsion. This result was improved by Wajnryb who found
the smallest possible generating set of $\Mod_{g,0}$ consisting of
$2$ elements one of which is torsion see~\cite {Wa}.
It is also known that in the case of closed surface the mapping
class group can be generated by 3 torsion elements (see~\cite{BF}),
two of which are involutions.
Recently Korkmaz (see~\cite{Kor}) showed that the mapping class
group can be generated by $2$ torsion elements (also in the case of
a closed surface).

In \cite{Mac}, Maclachlan proved that the moduli space ${\cal
M}_g$ is simply connected as a topological space by showing that
$\Mod_{g,0}$ is generated by torsion elements.  Several years
later Patterson generalized these results to $\Mod_{g,b}$ for
$g\geq 3,b\geq 1$ see \cite{Pa}. The question of generating
mapping class groups by involutions was considered by McCarthy and
Papadopoulos in \cite{MP}. Among other results, they proved that
for $g \geq 3$, $\Mod_{g,0}$ is generated by infinitely many
conjugates of a single involution.

Several year ago Luo, see \cite{Luo}, described the first finite
set of involutions which generate $\Mod_{g,b}$ for $g \geq 3$. The
size of his generating set depends linearly on both $g$ and $b$.
Luo also proved that $\Mod_{g,b}$ is generated by torsion elements
in all cases except $g=2$ and $b=5k+4$, but this group is not
generated by involutions if $g\leq 2$. In that paper Luo poses the
question of whether there is a universal upper bound, independent
of $g$ and $b$, for the number of torsion elements/involutions
needed to generate $\Mod_{g,b}$.

In the case of zero or one punctures, this question was answered by
Brendle and Farb in \cite{BF}, where they constructed a generating
set of $\Mod_{g,b}$ for $g\geq 3$ and $b\leq 1$ consisting of $7$
involutions. A detailed consideration of the generating set of
involutions constructed by Brendle and Farb shows that one of the
involutions is redundant and the mapping class group can be
generated by $6$ involutions (at least in the case $g>3$, or $g=3$
and no punctures).

Our main result generalizes the construction in the above paper to
the case of several punctures. We improve the bound given by Brendle
and Farb for the smallest number of involutions needed to generate
the mapping class group $\Mod_{g,b}$.

\begin{theorem}
\label{main} For all $g \geq 3$, the mapping class group
$\Mod_{g,b}$ can be generated by:
a) $4$ involutions if $g>7$ or $g=7$ and the number of punctures is even;

b) $5$ involutions if $g>5$ or $g=5$ and the number of punctures is even;

c) $6$ involutions if $g>3$ or $g=3$ and the number of punctures is even;

d) $9$ involutions if $g=3$ and the number of punctures is odd.
\end{theorem}

\begin{remark}
From the work of Luo \cite{Luo}, it follows that the above theorem can not be
generalized to the case of small genus, because the mapping class
group is not generated by involutions if $g \leq 2$.
\end{remark}

\begin{remark}
If one allows orientation-reversing involutions, it is possible to
extend Theorem \ref{main} to the extended mapping class group
which includes orientation-reversing mapping classes and show that
it is generated by the same number of involutions as the mapping
class group.
\end{remark}

\begin{remark}
Our result for genus $3$ is weaker than the result by Brendle and
Farb, but the construction in their paper does not work for $g=3$
and one puncture. In \cite{BF}, they wrote:
\begin{quote}
For simplicity of exposition, we provide explicit
arguments only for $b=0$.  In the case $b=1$ the arguments are the same,
although some involutions must be replaced with certain conjugates which
move the puncture to a fixed point of the involution.
\end{quote}
The problem with this argument is that for $g=3$ two of the
involutions in their generating set (the pair swap involutions
$J_1$ and $J_2$) do not have any fixed points.
\end{remark}

\begin{remark}
It is interesting to find the smallest possible set of involutions
which generated $\Mod_{g,b}$. It is clear that such a set should
contain at least $3$ involutions, since every group generated by
two involutions is dihedral.
One should expect that the bound
of $9$ involution for the genus $3$ case can be improved. The
author believes that the bound of $4$ involutions in the high
genus case could not be improved.
\end{remark}

{\em Acknowledgements:} The author wishes to thanks Prof. Sehgal, for
the useful discussions and the Department of Mathematics at
University of Alberta for the pleasant working environment.
\end{section}

\begin{figure}
\setlength{\unitlength}{0.01in}
\begin{picture}(0,0)(0,0)
\put(10,273){${}_{\alpha_1}$}
\put(77,238){${}_{\beta_1}$}
\put(35,280){${}_{\gamma_1}$}

\put(10,318){${}_{\alpha_2}$}
\put(77,320){${}_{\beta_2}$}

\put(170,273){${}_{\alpha_g}$}
\put(103,238){${}_{\beta_g}$}
\put(120,280){${}_{\gamma_{g-1}}$}

\put(165,320){${}_{\alpha_{g-1}}$}
\put(103,320){${}_{\beta_{g-1}}$}
\put(10,453){${}_{\alpha_{k-1}}$}
\put(67,410){${}_{\beta_{k-1}}$}
\put(25,460){${}_{\gamma_{k-1}}$}

\put(10,498){${}_{\alpha_k}$}
\put(57,503){${}_{\beta_k}$}

\put(160,457){${}_{\alpha_{k+2}}$}
\put(113,410){${}_{\beta_{k+2}}$}
\put(120,460){${}_{\gamma_{k+1}}$}

\put(160,498){${}_{\alpha_{k+1}}$}
\put(113,505){${}_{\beta_{k+1}}$}

\put(85,505){${}_{\gamma_{k}}$}
\put(-10,183){${}_{x_1}$}
\put(-10,150){${}_{x_2}$}

\put(190,183){${}_{x_b}$}
\put(190,150){${}_{x_{b-1}}$}

\put(-20,70){${}_{x_{l-1}}$}
\put(190,70){${}_{x_{l+1}}$}
\put(80,40){${}_{x_{l}}$}
\put(85,185){${}_{\delta}$}
\put(90,560){${}_{\rho_1}$}
\put(290,253){${}_{\alpha_1}$}
\put(375,245){${}_{\beta_1}$}
\put(353,297){${}_{\gamma_1}$}

\put(270,318){${}_{\alpha_2}$}
\put(330,320){${}_{\beta_2}$}

\put(455,318){${}_{\alpha_g}$}
\put(395,320){${}_{\beta_g}$}
\put(270,428){${}_{\alpha_{k-1}}$}
\put(327,388){${}_{\beta_{k-1}}$}
\put(290,435){${}_{\gamma_{k-1}}$}

\put(270,473){${}_{\alpha_k}$}
\put(295,483){${}_{\beta_k}$}
\put(318,485){${}_{\gamma_k}$}

\put(445,428){${}_{\alpha_{k+3}}$}
\put(388,388){${}_{\beta_{k+3}}$}
\put(399,435){${}_{\gamma_{k+2}}$}

\put(445,473){${}_{\alpha_{k+2}}$}
\put(370,448){${}_{\beta_{k+2}}$}
\put(407,484){${}_{\gamma_{k+1}}$}

\put(385,508){${}_{\beta_{k+1}}$}
\put(325,508){${}_{\alpha_{k+1}}$}
\put(253,183){${}_{x_1}$}
\put(253,150){${}_{x_2}$}

\put(475,183){${}_{x_{b-1}}$}
\put(475,150){${}_{x_{b-2}}$}

\put(243,70){${}_{x_{l-1}}$}
\put(475,70){${}_{x_{l}}$}

\put(370,270){${}_{x_{b}}$}
\put(360,185){${}_{\delta}$}
\put(365,560){${}_{\rho_2}$}
\end{picture}
\includegraphics{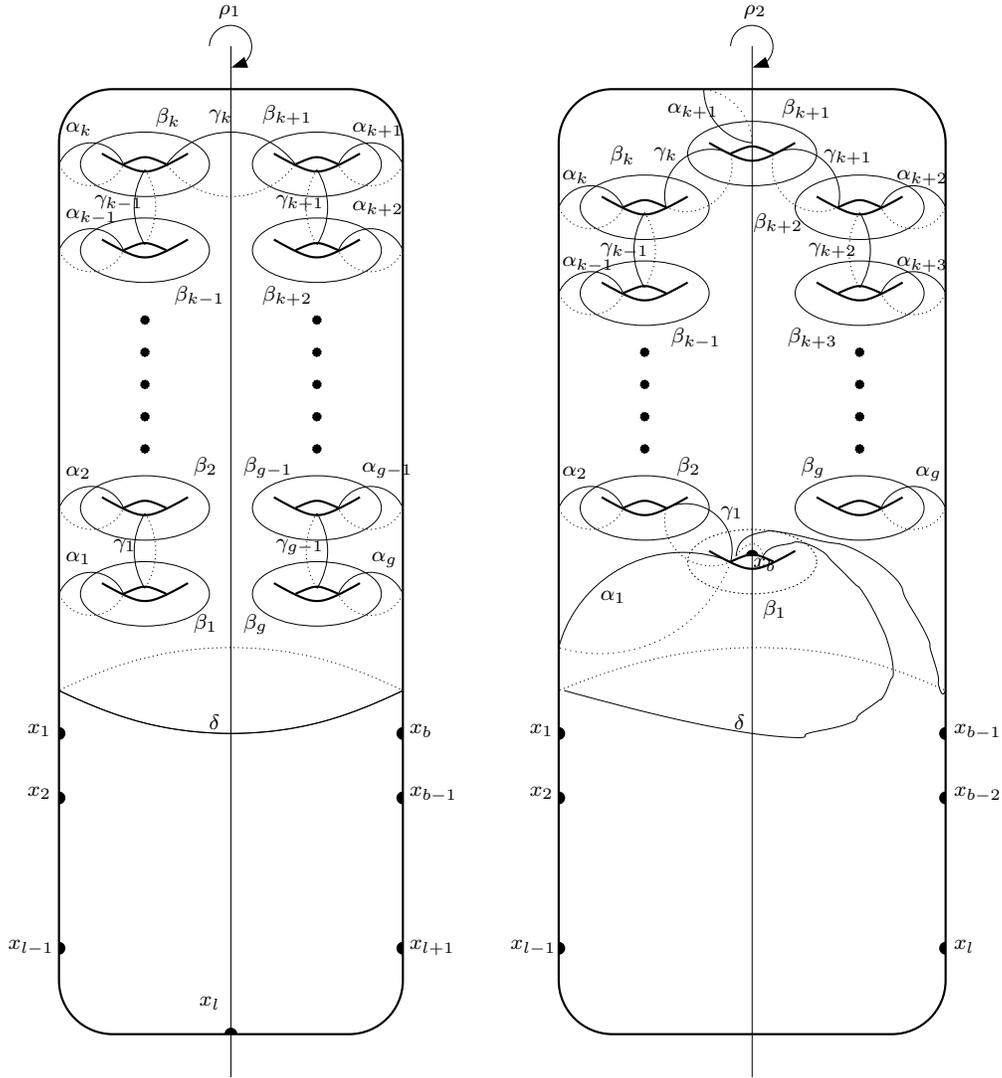}%
\caption{The embeddings of the surface $\Sigma_{g,b}$ in the Euclidian
space used to define the involutions $\rho_1$ and $\rho_2$.}
\label{2embedings}
\end{figure}

\begin{section}{Generating the mapping class group by 2
involutions and 3 Dehn twists}

Let us embed our surface $\Sigma=\Sigma_{g,b}$ in the Euclidian space
in two different ways as shown on Figure~\ref{2embedings}.
(In these pictures we will assume that genus $g=2k$ is even and the
number of punctures $b=2l+1$ is odd. In the case of odd genus we only
have to swap the top parts of the pictures, and in the case of
even number of punctures we have to swap the lower parts of the
pictures.)

In Figure~\ref{2embedings} we have also marked the puncture points
as $x_1, \dots, x_b$ and we have the curves $\alpha_i$,
$\beta_i$, $\gamma_i$ and $\delta$. Using these curves one can
construct an explicit homeomorphism between the two surfaces in
Euclidian space. The curve $\delta$ separates our surface into two
components: the first one, denoted by $\Sigma'$, is a surface of
genus $g$ with boundary component and no punctures. The second one
(denoted by $D$) is a disk with $b$ puncture points. All other
curves are simple non-separating curves on the surface $\Sigma'$.

Each embedding gives a natural involution of the surface --- the half turn
rotation around its axis of symmetry. Let us call these involutions
$\rho_1$ and $\rho_2$. The product $R=\rho_2\rho_1$ of these involutions acts
almost as a rotation on the surface $\Sigma_{g,b}$\footnote{In the case of
no punctures $R$ is homeomorphic to a rotation if one uses the flower
embedding of the surface $\Sigma_{g,0}$ in the Euclidian space}, more
precisely we have
\begin{equation}
\label{rotation}
\begin{array}{cl}
R \alpha_i = \alpha_{i+1}, & \mbox{ for } 1 \leq i < g \\
R \beta_i  = \beta_{i+1},  & \mbox{ for } 1 \leq i < g \\
R \gamma_i = \gamma_{i+1}, & \mbox{ for } 1 \leq i < g-1\footnotemark.
\end{array}
\end{equation}
\footnotetext{We also have that $R\beta_g=\beta_1$, but
$R\alpha_g\not=\alpha_1$ unless there are no punctures.}
On the set of punctures $R$ acts also as a long cycle
$$
Rx_1 = x_b \quad \mbox{and} \quad Rx_i = x_{i-1} \mbox{ for } 1 < i \leq b.
$$

Let $c$ be a simple non-separating curve on an oriented surface
$\Sigma$.  We will denote by $T_c$ the Dehn twist around the curve
$c$\footnote{We will not distinguish between the homeomorphisms
and their images in the mapping class group $\Mod_{g,b}$.}.
For any homeomorphism $h$ of the surface $\Sigma$ the twists
around the curves $c$ and
$h(c)$ are conjugate in the mapping class group $\Mod(\Sigma)$,
$$
T_{h(c)} = h T_c h^{-1}.
$$

The relations (\ref{rotation}) imply that the Dehn twists around the
curves $\alpha_i$-es, $\beta_i$-es and $\gamma_i$-es are
conjugate, i.e.,
\begin{equation}
\label{conjugation}
R T_{\alpha_i} R^{-1} = T_{\alpha_{i+1}}, \quad
R T_{\beta_i}  R^{-1} = T_{\beta_{i+1}} \mbox{ and }
R T_{\gamma_i} R^{-1} = T_{\gamma_{i+1}}.
\end{equation}

The Dehn twists around these $3g-1$ curves generate the
whole mapping class group of a closed surface of genus $g$ by a result of
Lickorish, and  of the surface $\Sigma'$ of genus $g$ and one boundary component
by a result of Johnson (see~\cite{Li,Jo}).
This together with relations (\ref{conjugation}),
shows the the subgroup $G$ generated by
$\rho_1$, $\rho_2$ and $3$ Dehn twists $T_\alpha$, $T_\beta$ and $T_\gamma$
around one of the curve in each family, contains
$$
\Mod(\Sigma') \subset \Mod(\Sigma_{g,b}).
$$

Our next step is to show that this subgroup actually contains all
homeomorphisms of the surface $\Sigma_{g,b}$ which preserve the
punctures pointwise. In~\cite{Ge} it is shown that the group
$\Mod^0(\Sigma_{g,b})$ is generated by the Dehn twists around the
curves $\alpha_i$-es, $\beta_i$-es, $\gamma_i$-es and
$\delta_j$-es, for $j=1,\dots,b-1$, where the curves $\delta_i$-es are
shown on Figure~\ref{deltai-es}.

\begin{figure}
\setlength{\unitlength}{0.01in}
\begin{picture}(0,0)(0,0)
\put(10,273){${}_{\alpha_1}$}
\put(170,273){${}_{\alpha_g}$}

\put(50,283){${}_{\beta_1}$}
\put(130,283){${}_{\beta_g}$}
\put(-10,183){${}_{x_1}$}
\put(-10,150){${}_{x_2}$}

\put(190,183){${}_{x_b}$}
\put(190,150){${}_{x_{b-1}}$}

\put(-20,70){${}_{x_{l-1}}$}
\put(190,70){${}_{x_{l+1}}$}
\put(80,40){${}_{x_{l}}$}

\put(80,450){${}_{\rho_{1}}$}
\put(-10,210){${}_{\delta_0}$}
\put(-10,170){${}_{\delta_1}$}

\put(190,205){${}_{\delta_b}$}
\put(190,165){${}_{\delta_{b-1}}$}

\put(65,5){${}_{\delta_{l-1}}$}
\put(115,5){${}_{\delta_{l}}$}
\put(260,273){${}_{\alpha_1}$}
\put(420,273){${}_{\alpha_g}$}

\put(300,283){${}_{\beta_1}$}
\put(380,283){${}_{\beta_g}$}
\put(240,183){${}_{x_1}$}
\put(240,150){${}_{x_2}$}

\put(440,183){${}_{x_b}$}
\put(440,150){${}_{x_{b-1}}$}

\put(230,70){${}_{x_{l-1}}$}
\put(440,70){${}_{x_{l+1}}$}
\put(330,40){${}_{x_{l}}$}
\put(240,210){${}_{\eta_0}$}
\put(240,170){${}_{\eta_1}$}

\put(440,205){${}_{\eta_b}$}
\put(440,165){${}_{\eta_{b-1}}$}

\put(315,5){${}_{\eta_{l-1}}$}
\put(365,5){${}_{\eta_{l}}$}

\put(330,450){${}_{\rho_{1}}$}
\end{picture}
\includegraphics{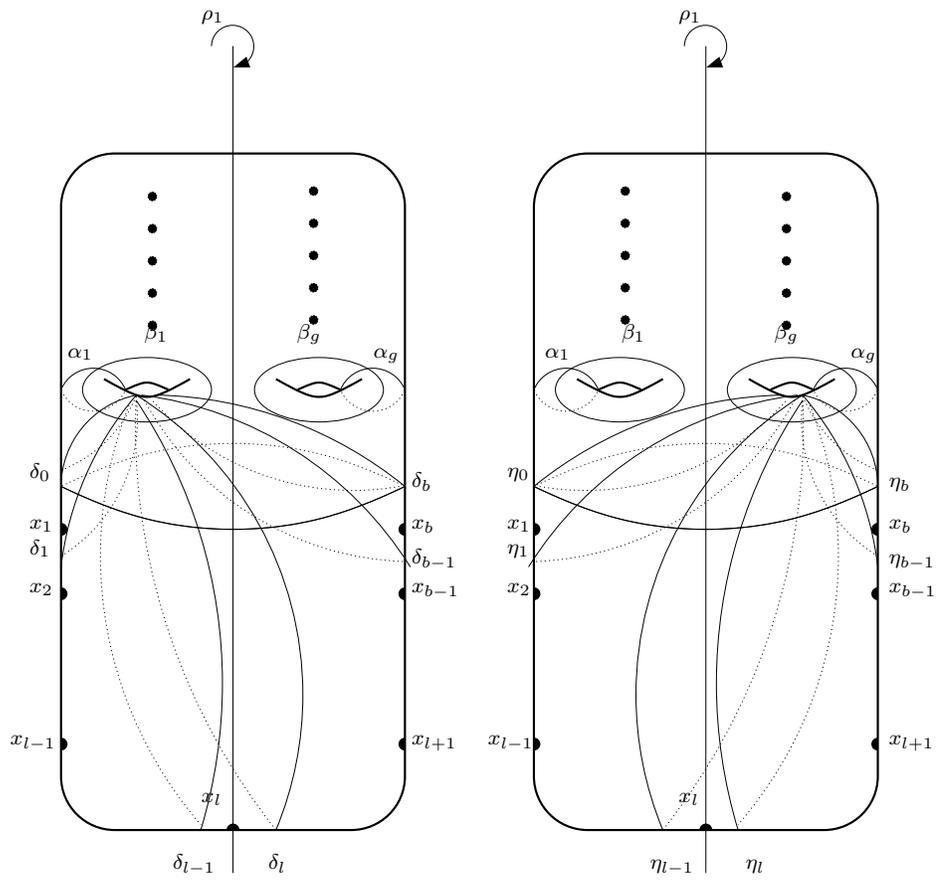}%
\caption{The curves $\delta_i$-es and $\eta_i$.}
\label{deltai-es}
\end{figure}

\begin{lemma}
\label{Rdeltaeta}
The `rotation' $R$ acts on the curves $\delta_j$ as follows:
$$
R^{-1} \delta_i = \eta_{i+1}
$$
\end{lemma}
\begin{proof}
The Figure \ref{Rondeltai-es} shows the action of $\rho_1$ and $\rho_2$
on the curve $\delta_i$. It is clear from the picture that
$\eta_{i+1}= \rho_1 \rho_2 (\delta_i) = R^{-1}\delta_i$.
\begin{figure}
\setlength{\unitlength}{0.01in}
\begin{picture}(0,0)(0,0)
\put(-5,0){
\begin{picture}(0,0)(0,0)
\put(10,93){${}_{x_i}$}
\put(10,73){${}_{x_{i+1}}$}
\put(10,53){${}_{x_{i+2}}$}
\put(70,110){${}_{\delta}$}
\put(30,100){${}_{\delta_i}$}
\put(10,160){${}_{\alpha_1}$}
\put(35,140){${}_{\beta_1}$}
\put(80,160){${}_{\alpha_g}$}
\put(50,140){${}_{\beta_g}$}
\put(40,240){${}_{\rho_1}$}
\end{picture}
}
\put(45,0){$\delta_i$}
\put(95,130){$=$}
\put(100,0){
\begin{picture}(0,0)(0,0)
\put(10,93){${}_{x_i}$}
\put(10,73){${}_{x_{i+1}}$}
\put(10,53){${}_{x_{i+2}}$}
\put(70,110){${}_{\delta}$}
\put(35,100){${}_{\delta_i}$}
\put(10,160){${}_{\alpha_1}$}
\put(55,140){${}_{\beta_1}$}
\put(40,240){${}_{\rho_2}$}
\end{picture}
}
\put(150,0){$\delta_i$}
\put(205,0){
\begin{picture}(0,0)(0,0)
\put(65,93){${}_{x_{b-i}}$}
\put(55,73){${}_{x_{b-i-1}}$}
\put(55,53){${}_{x_{b-i-2}}$}
\put(20,105){${}_{\delta}$}
\put(50,120){${}_{\rho_2(\delta_i)}$}
\put(10,160){${}_{\alpha_1}$}
\put(55,140){${}_{\beta_1}$}
\put(40,240){${}_{\rho_2}$}
\end{picture}
}
\put(250,0){$\rho_2(\delta_i)$}
\put(305,130){$=$}
\put(310,0){
\begin{picture}(0,0)(0,0)
\put(53,93){${}_{x_{b-i+1}}$}
\put(65,73){${}_{x_{b-i}}$}
\put(53,53){${}_{x_{b-i-1}}$}
\put(75,115){${}_{\delta}$}
\put(20,80){${}_{\rho_2(\delta_i)}$}
\put(10,160){${}_{\alpha_1}$}
\put(35,140){${}_{\beta_1}$}
\put(80,160){${}_{\alpha_g}$}
\put(50,140){${}_{\beta_g}$}
\put(40,240){${}_{\rho_1}$}
\end{picture}
}
\put(350,0){$\rho_2(\delta_i)$}
\put(415,0){
\begin{picture}(0,0)(0,0)
\put(10,93){${}_{x_i}$}
\put(10,73){${}_{x_{i+1}}$}
\put(10,53){${}_{x_{i+2}}$}
\put(70,110){${}_{\delta}$}
\put(40,100){${}_{\eta_i}$}
\put(10,160){${}_{\alpha_1}$}
\put(35,140){${}_{\beta_1}$}
\put(80,160){${}_{\alpha_g}$}
\put(50,140){${}_{\beta_g}$}
\put(40,240){${}_{\rho_1}$}
\end{picture}
}
\put(415,0){$\rho_1\rho_2(\delta_i)=\eta_{i+1}$}
\end{picture}
\includegraphics{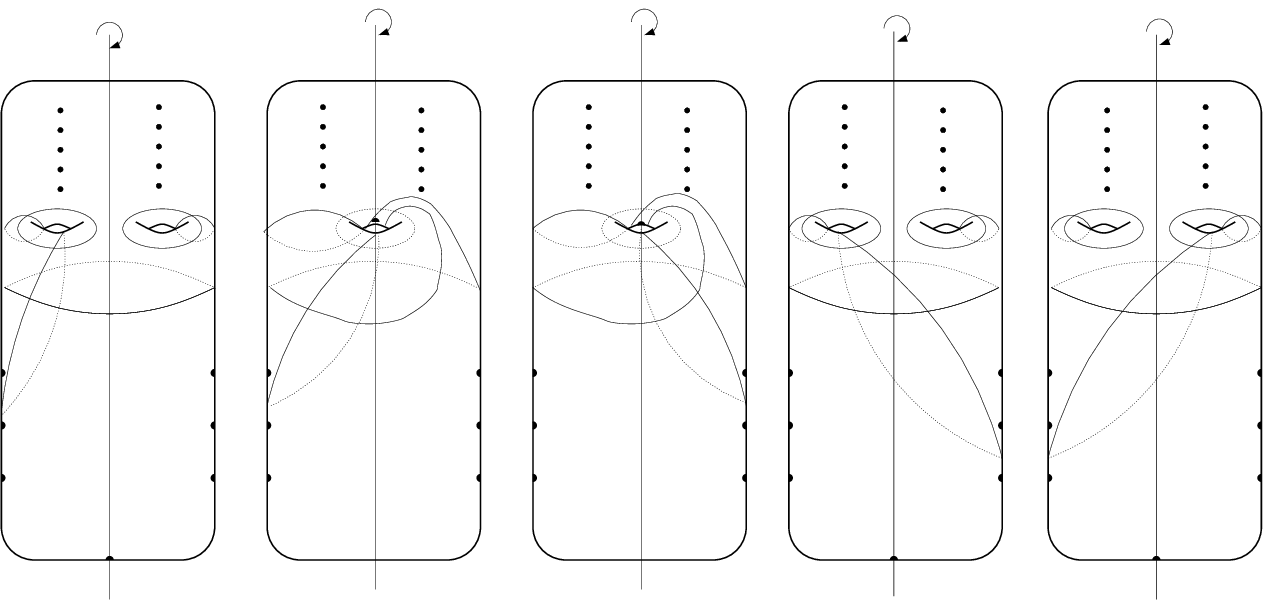}%
\caption{The action of $R$ on the curve $\delta_i$.}
\label{Rondeltai-es}
\end{figure}
\end{proof}

Now we are ready to show that $G$ contains $\Mod^0(\Sigma_{g,b})$.

\begin{lemma}
The group $G$ contains the Dehn twists around the curves
$\delta_j$, for $j=0,\dots, b-1$.
\end{lemma}
\begin{proof}
We will prove the lemma using induction on $j$.
The base case, $j=0$, is clear because by construction
$G$ contains $T_{\delta_0} = T_{\alpha_1} \in G$. Suppose that $G$
contains the twist $T_{\delta_j}$.
Using Lemma \ref{Rdeltaeta} we can see that the
twist $T_{\eta_{j+1}}$ also lies in $G$, since it is
conjugate to $T_{\delta_j}$
$$
T_{\eta_{j+1}} = R^{-1} T_{\delta_j} R \in G.
$$
To complete the induction step, we only need to notice that there exists
a homeomorphism $h_{j+1} \in \Mod(\Sigma') \subset G $ such that
$h_{j+1} \eta_{j+1} = \delta_{j+1}$, because the curves
$\delta_{j+1}$ and $\eta_{j+1}$ are homotopic in the disk $D$.
This gives us that
$$
T_{\delta_{j+1}} = h_{j+1}T_{\eta_{j+1}} h_{j+1}^{-1} \in G,
$$
which completes the induction step.
\end{proof}

\begin{corollary}
\label{mod0}
The group $G$ contains the subgroup $\Mod^0(\Sigma_{g,b})$.
\end{corollary}

\begin{remark}
The group $G$ is not the whole group $\Mod^0(\Sigma_{g,b})$ (if
$b>3$), because its image in the symmetric group $\Sym_b$ is
generated by two involutions and therefore is a dihedral group. It
is easy see that this image is the group $D_{2b}$ which is a proper subgroup
of $\Sym_b$. However in the next section we will show how to
generated the Dehn twists $T_\alpha$, $T_\beta$ and $T_\gamma$ by
several involutions and we can be flexible in choosing their
actions on the set of punctures. We can choose these involutions so
that their images in the symmetric group together with the images
of $\rho_1$ and $\rho_2$ generate the whole symmetric group on the
punctures.
\end{remark}
\end{section}

\begin{section}{Generating Mapping class group by involutions}

In the previous section we have shown that we can generate almost the whole
mapping class group by $2$ involutions and $3$ Dehn twists. Now we only need
to explain how to generate the Dehn twists by involutions. The basic idea is
to use the lantern relation and generate the Dehn twist $T_\alpha$ by
four involutions.

\begin{subsection}{The Lantern Relation}

We begin by recalling the lantern relation in the mapping class
group. This relation was first discovered by Dehn in the 1930s and
later rediscovered by Johnson. Here we will use this relation to
show how a Dehn twist can be expressed using $4$ involutions. This
argument is based on a result by Luo and Harer (see~\cite{Luo,Ha}), who
proved that a Dehn twist lie in a subgroup generated by $6$ involutions.
Brendle and Farb in~\cite{BF} improved this argument by
introducing pair swap involutions and showed that $4$ involutions
suffice to generate a Dehn twist.

Let the $S_{0,4}$ be a surface of genus $0$ with $4$ boundary
components. Denote by $a_1, a_2, a_3$ and $a_4$ the four boundary
curves of the surface $S_{0,4}$ and let the interior curves $x_1,
x_2$ and $x_3$, be as shown in Figure \ref{lantern}.

\begin{figure}
\center{
\setlength{\unitlength}{0.01in}
\begin{picture}(0,0)(0,0)
\put(0,0){
\begin{picture}(0,0)(0,0)
\put(80,150){${}_{a_1}$}
\put(80,10){${}_{a_2}$}
\put(150,80){${}_{a_3}$}
\put(5,80){${}_{a_4}$}
\put(35,50){${}_{x_1}$}
\put(35,125){${}_{x_2}$}
\put(135,55){${}_{x_3}$}
\end{picture}
}

\put(250,80){$J_{12}$}

\put(305,0){
\begin{picture}(0,0)(0,0)
\put(80,150){${}_{a_1}$}
\put(80,10){${}_{a_2}$}
\put(150,80){${}_{a_3}$}
\put(5,80){${}_{a_4}$}
\put(35,50){${}_{x_1}$}
\put(35,125){${}_{x_2}$}
\put(135,55){${}_{x_3}$}
\end{picture}
}

\end{picture}
\includegraphics{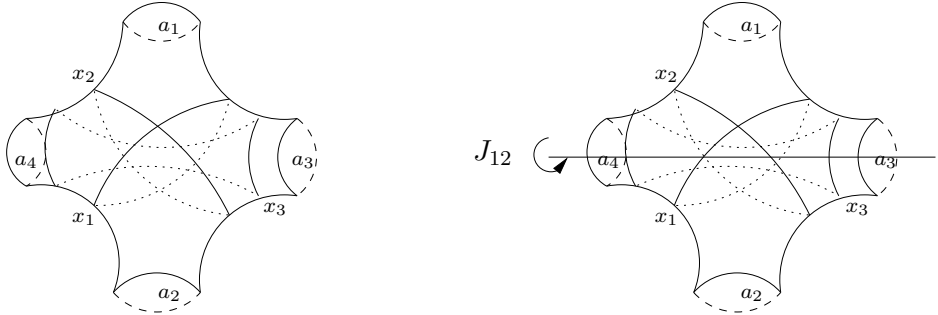}
\caption{Lantern and a pair swap involution}
}
\label{lantern}
\label{pairswap}
\end{figure}

The following relation :
\begin{equation}
\label{eq:lantern}
T_{x_1} T_{x_2}T_{x_3} = T_{a_1} T_{a_2} T_{a_3} T_{a_4}
\end{equation}
among the Dehn twists around the curves $a_i$ and $x_i$ is known
as the lantern relation. Notice that the curves $a_i$ do not
intersect any other curve and that the Dehn twists $T_{a_i}$
commute with every twist in this relation. This allows us to
rewrite the lantern relation as follows
\begin{equation}
\label{eq:lantern1}
T_{a_4} = (T_{x_1} T_{a_1}^{-1}) (T_{x_2}T_{a_2}^{-1}) (T_{x_3}T_{a_3}^{-1})
\end{equation}

The product $T_{x_1} T_{a_1}^{-1}$ can be expressed as a product
of two involutions in the following way: let $\rho$ be an
involution which maps $a_1$ to $x_1$, then we have
$$
T_{x_1} = \rho T_{a_1} \rho
$$
and
\begin{equation}
\label{lantern2}
T_{x_1} T_{a_1}^{-1} = \rho T_{a_1} \rho T_{a_1}^{-1}
= \rho ( T_{a_1} \rho T_{a_1}^{-1})
\end{equation}

Figure \ref{pairswap} shows the  existence of a pair swap
involutions $J_{12}$ on $\Sigma_{0,4}$, such that
$J_{12}$ preserves the curves $a_4$ and $a_3$, and swaps the pairs
$(a_1, x_1)$ with $(a_2,x_2)$. There also exists another pair
swap involution $J_{13}$ with similar properties.
This gives that
$$
\begin{array}{c}
T_{x_2} T_{a_2}^{-1} = J_{12} T_{x_1} T_{a_1}^{-1} J_{12} \\
T_{x_3} T_{a_3}^{-1} = J_{13} T_{x_1} T_{a_1}^{-1} J_{13}
\end{array}
$$

These relations together with relation \ref{eq:lantern1} and
\ref{lantern2} show that the twist around the curve $a_4$ lies in
the group generated by the involutions $\rho$, $\bar \rho =T_{a_1}
\rho T_{a_1}^{-1}$, $J_{12}$ and $J_{13}$, more precisely we have
$$
T_{a_4} =
(\rho \bar \rho )
(J_{12} \rho \bar \rho J_{12})
(J_{13} \rho \bar \rho J_{13}).
$$
\end{subsection}

After this preparation, we can start the prof of Theorem~\ref{main}.
We will begin with part c).

\begin{subsection}{Genus $g\geq 3$ and $b$ even if $g=3$}

The argument which follows is the same as the one used in \cite{BF}.
The only major difference is that we use a different lantern.

Recall that in the whole paper we assume that the genus of the
surface is at least $3$. In the argument which follows we will
need an additional assumption that the number of punctures is
even, if the genus $g$ is exactly $3$. For simplicity we will assume
that the genus $g=2k$ is even. In the case of odd genus $g=2k+1$,
we can use the same lantern but use the involution $\rho_2$ instead of
$\rho_1$.

We want to apply the argument in the previous section and find $4$
involutions which generate the twist $T_\alpha$.
In order to do that we need to find an embedding of
the lantern $S_{0,4}$ into $\Sigma_{g,b}$
and extend all involutions to the whole surface $\Sigma_{g,b}$.

\begin{figure}
\center{
\setlength{\unitlength}{0.01in}
\begin{picture}(0,0)(0,0)

\put(180,255){$_{\alpha_{k+1}}$}
\put(20,255){$_{\alpha_{k}}$}
\put(15,215){$_{\alpha_{k-1}}$}
\put(95,260){$_{\gamma_k}$}
\put(60,240){$_{\gamma_{k-1}}$}
\put(100,340){$_{\rho_1}$}

\put(300,255){$_{\alpha_k}$}
\put(440,260){$_{\beta_{k+1}}$}
\put(370,270){$_{x}$}
\put(380,340){$_{\rho_1}$}
\end{picture}
\includegraphics{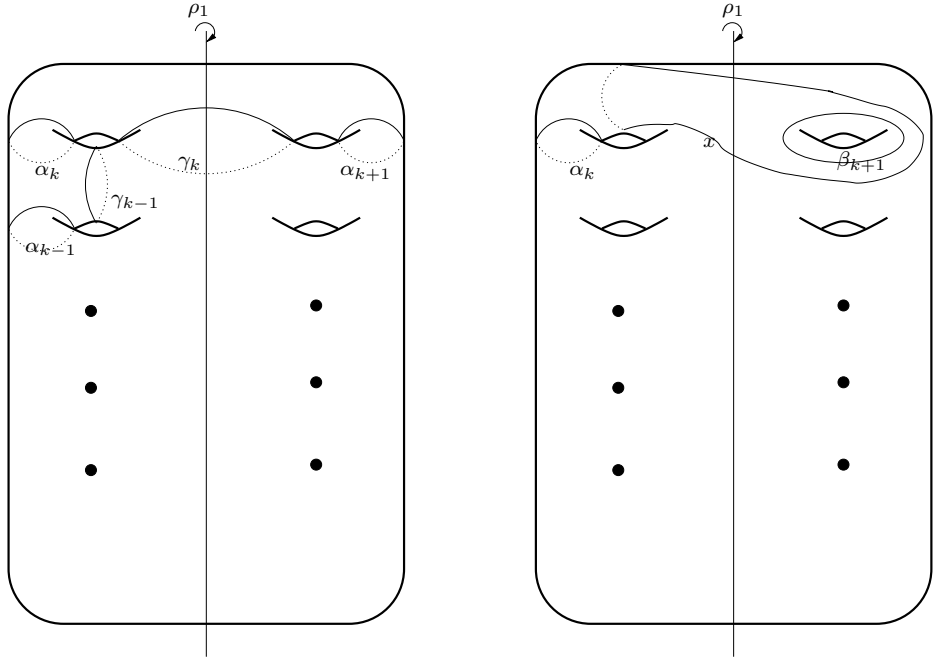}
}
\caption{Embedding of a lantern and a pair of pants in $\Sigma_{g,b}$.}
\label{lanterninsigma}
\end{figure}

The Figure~\ref{lanterninsigma} shows one embedding of a lantern
in our surface $\Sigma_{g,b}$. Note that the compliment of this
lantern is still a connected surface of genus $g-3$ with $4$
boundary components and $b$ punctures. The boundary components of
this lantern are $a_1 = \alpha_{k+1}$, $a_2 = \alpha_{k-1}$, $a_3 =
\gamma_{k-1}$ and $a_4 = \gamma_{k}$; and the middle curve $x_1$ is
$\alpha_{k}$. Let denote $S$ denotes this embedding of a lantern in $\Sigma$.

Let us notice that the involution $\rho_1$ ($\rho_2$ in the case of odd $g$)
 maps $a_1$ to $x_1$. Let $\rho_3$ denotes the
 conjugate of the involution $\rho_1$ (or $\rho_2$)
with $T_{x_1}$. This gives us two of the four
involutions needed to express $T_{a_4} = T_{\gamma_{k}}$.

The compliment of the lantern is a surface $S_1$ of genus $g-3$ with
$4$ boundary components and $b$ punctures.

On a surface $S_1'$ (of genus $g-3$ with $4$ boundary components)
there exists an involution $\tilde J_{12}$ which preserves the two
boundary components $a_3$ and $a_4$ and switches the other two.
Moreover, we can arrange that this involution has
$2(g-3)$ fixed points on the surface $S_1'$. We can find an invariant
set (under the action of the involution $\tilde J_{12}$) of
$b$ points/punctures $S_1'$, here we use
the condition that $g>3$ if $b$ is odd, which implies that
there exist a finite $\tilde J_{12}$ invariant set with arbitrary
number of points. This allows us to define
an involution $\bar J_{12}$ on the surface $S_1$, which preserves
two of the boundary components and switches the other two.
Finally, we can define an involution $I_{12}$ on the surface
$\Sigma_{g,b}$ by gluing together the involutions $J_{12}$ and
$\bar J_{12}$ on the two parts of surface. In exactly the same way
we can construct an involution $I_{13}$.

From the above construction it is clear that the twist
$T_{\gamma_{k}}$ lies in the group $G_1$ generated by the
involutions $\rho_1$, $\rho_2$, $\rho_3$, $I_{12}$ and $I_{13}$.
Using conjugation with the element $R$ we can see that $G_1$
contains the twists around all curves $\gamma_i$-es.

By construction the involution $I_{13}$
takes the curve $\alpha_{k+1}$ to the curve $\gamma_{k-1}$, i.e., we
have
$$
T_{\alpha_{k+1}} = I_{13} T_{\gamma_{k-1}} I_{13} =
I_{13} R T_{\gamma_{k}} R^{-1} I_{13} \in G_1.
$$

So far we have used $5$ involutions to generate a group $G_1$
which contains the twists around the curves $\alpha_i$-es and
$\gamma_i$-es.

As shown on Figure \ref{lanterninsigma},
we can find a pair of pants (surface of genus $0$,
with $3$ boundary components) on the surface $\Sigma_{g,b}$, with
boundary curves $\alpha_{k}$, $\beta_{k+1}$ and some curve $x$.

Using the same construction as above, we can find an involution
$I$ on the surface $\Sigma_{g,b}$ which interchanges the curves
$\alpha_k$ and $\beta_{k-1}$ and fixes the curve $x$.
Moreover, since the compliment of a
pair of pants is a surface of genus $g-2 \geq 1$, we can chose $I$
so that it acts on the set of punctures as any involution in $\Sym_b$ with up
to $3$ fixed points.

\begin{lemma}
The group $G_2$ generated by the involutions $\rho_1$, $\rho_2$, $\rho_3$,
$I_{12}$, $I_{13}$ and $I$ contains $\Mod^0(\Sigma_{g,b})$.
\end{lemma}
\begin{proof}
We have already shown that $G_1 \subset G_2$ contains the Dehn twists
around the curves $\alpha_i$-es and $\gamma_i$-es. Using the conjugation
by the involution $I$, we can see that $G_2$ also contains the twist around
$\beta_{i-1}$. This allows us to apply corollary~\ref{mod0} and conclude that
$G_2$ contains the subgroup $\Mod^0(\Sigma{g,b})$.
\end{proof}

Knowing that $G_2$ contains $\Mod^0(\Sigma_{g,b})$, we only need
to make sure that all permutations of the punctures can be
obtained
as elements of $G_2$ to conclude that $G_2$ is the full mapping
class group.

\begin{theorem}[6 involutions]
\label{6inv}
If $g\geq 3$, and $b$ is even if $g=3$, we can choose
the involution $I$ such that the group generated by $\rho_1$,
$\rho_2$, $\rho_3$, $I_{12}$, $I_{13}$ and $I$ is the whole
mapping class group $\Mod_{g,b}$.
\end{theorem}
\begin{proof}
To prove the theorem we only need to show that $G_2$ can be made to map
surjectively onto the symmetric group $\Sym_b$,
because $G_2$ contains $\Mod^0(\Sigma_{g,b})$.
This is equivalent to showing that the full symmetric group $\Sym_b$
can be generated by the involutions
$r_1=(1,b)(2,b-1)(3,b-2)\dots$,
$r_2=(1,b-1)(2,b-2)\dots$ corresponding to $\rho_1$ and $\rho_2$ and another
involution which have less then $4$ fixed points.

\begin{lemma}
\label{genofSb}
The symmetric group $\Sym_b$ is generated by $r_1$, $r_2$ and the
involution $r_3=(2,b-1)(3,b-2)\dots$.
\end{lemma}
\begin{proof}
The group generated by $r_i$ contains the long cycle $r_1r_2=(12\dots b)$
and the transposition $r_1r_3=(1,b)$. These two elements generate the
whole symmetric group, therefore the involutions $r_i$ generate $\Sym_b$.
\end{proof}

Thus we have shown that if $I$ acts on the punctures as $r_3$ then the
group generated by $\rho_1$, $\rho_2$, $\rho_3$,
$I_{12}$, $I_{13}$ and $I$ is $\Mod_{g,b}$, which finishes
the proof of part c) of Theorem~\ref{main}.
\end{proof}
\end{subsection}

\begin{subsection}{Genus at least 5}

Now we would like to refine the above argument and show that for
large genus, we can avoid using the involution $I$, in order to
generate the mapping class group. From now on we will assume that the
genus $g$ of the surface is at least $5$ (and that the number of punctures
is even in the case $g=5$).

We will start with the same lantern as above and define the
involution $I_{13}$ in the same way. Now we know that the genus
$g$ is at least $5$ which allows us to make the action of $I_{13}$
on the set of punctures the same as the action of any involution
with up to $3$ fixed points. In particular we can make $I_{13}$ to
act on the set of punctures as the involution $r_3$ used in the proof of
Lemma~\ref{genofSb}.

We will slightly modify the construction involution $I_{12}$.
The purpose of this modification is to make sure that $I_{12}$
sends some of the $\alpha_i$-es to some of the $\beta_i$-es, which
will allow us to remove the involution $I$ for the generating set.

Let us glue to our lantern, two pairs of pants, to obtain a
surface $S_2$ homeomorphic to a sphere with six boundary
components. The first pair of pants is bounded by $a_1 =
\alpha_{k+1}$, $\alpha_{k+2}$ and $\gamma_{k+1}$. The boundary components
of the  second pair of pants are $a_2 = \alpha_{k-1}$,
$\beta_{k-2}$ and some curve $x$ (see Figure \ref{2pairofpants}).
From the figure is clear that the complement of $S_2$ in
$\Sigma_{g,b}$ is a connected surface $S_3$ of genus $g-5$ with
$6$ boundary components and $b$ punctures.

\begin{figure}
\center{
\setlength{\unitlength}{0.0083in}
\begin{picture}(0,0)(0,0)
\put(100,240){${}_{\alpha_{k-1}}$}
\put(100,150){${}_{\alpha_{k+1}}$}
\put(45,210){${}_{\gamma_{k-1}}$}
\put(135,210){${}_{\gamma_{k}}$}

\put(135,280){${}_{\alpha_{k+2}}$}
\put(135,320){${}_{\gamma_{k+1}}$}

\put(135,80){${}_{\beta_{k-2}}$}
\put(135,120){${}_{x}$}

\put(0,215){${}_{J_{12}}$}

\put(230,305){${}_{\alpha_{k-1}}$}
\put(425,355){${}_{\alpha_{k+1}}$}
\put(290,310){${}_{\gamma_{k-1}}$}
\put(350,360){${}_{\gamma_{k}}$}

\put(425,305){${}_{\alpha_{k+2}}$}
\put(370,310){${}_{\gamma_{k+1}}$}

\put(275,250){${}_{\beta_{k-2}}$}
\put(320,200){${}_{x}$}

\put(335,415){${}_{\rho_1}$}

\end{picture}
\includegraphics{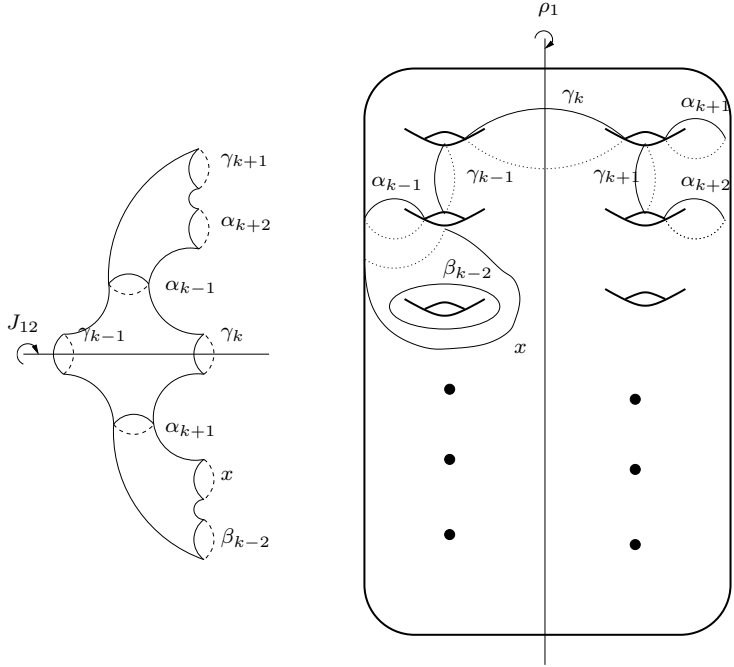}
}
\caption{Extending a lantern with two pair of pants.}
\label{2pairofpants}
\end{figure}

It is clear that we can extend the involution $J_{12}$ to a
involution $\hat J_{12}$ of the surface $S_2$ which takes
$\gamma_{k+1}$ to $\beta_{k-2}$ and $\alpha_{k+2}$ to $x$, and
preserves the other two boundary components,
$a_3=\gamma_{k-1}$ and $a_4=\gamma_{k}$. As before, we can
construct an involution $\bar J_{12}$ on the $S_3$ which preserves
two of the boundary components and swaps the other $4$ components
into two pairs. Gluing together $\hat J_{12}$ and $\bar J_{12}$
gives us the involution $I_{12}$ of the surface $\Sigma_{g,b}$.

\begin{theorem}[5 involutions]
\label{5inv}
If $g\geq 5$ and $b$ is even if $g=5$,
the group generated by $\rho_1$, $\rho_2$, $\rho_3$,
$I_{12}$ and $I_{13}$
is the whole mapping class group $\Mod_{g,b}$.
\end{theorem}
\begin{proof}
By construction the group $G_3$ generated by these $5$ involutions
contains the Dehn twist around the curve $\gamma_{k}$ and using
conjugation with powers of $R$ we have that $G_3$ contains the
twists around all $\gamma_i$-es. We have constructed the
involution $I_{12}$ in such a way that it sends one of the curves
$\gamma_i$-es into some $\beta_{k-2}$,
also the involution $I_{13}$ sends some $\gamma_i$ to $\alpha_{k+1}$.
This gives us that the Dehn twists
$T_{\beta_{k-2}}$ and $T_{\alpha_{k+1}}$ lie in the group $G_3$.
Using the Corollary~\ref{mod0} we can conclude that $G_3$ contains
the group $\Mod^0(\Sigma_{g,b})$.

We also have that the group $G_3$ maps onto $\Sym_b$ because its image
contains the involutions $r_1$, $r_2$ and $r_3$, which proves that
$G_3$ is the whole mapping class group.
\end{proof}
\end{subsection}

\begin{subsection}{High genus}
Now we want to refine the above argument and show that for high
genus we do not need the involution $I_{13}$ in order to generate
the mapping class group. Assume that the genus of the surface is
at least $7$ and that the number of puncture points is even if $g=7$.

We will construct an involution $I$, such that in the group generated by
$I$ and $R$ there are elements $h_{12}$ and $h_{13}$ such that
$$
h_{12}a_1 = a_2, \,\,\, h_{12}x_1 = x_2; \,\,\, \mbox{and} \,\,\,\\
h_{13}a_1 = a_3, \,\,\, h_{13}x_1 = x_3.
$$
This would imply that the group generated by $\rho_i$ and $I$
contains the product of twists
$T_{x_1}T^{-1}_{a_1} = \rho_1 \rho_3$ and its conjugates
$$
T_{x_2}T^{-1}_{a_2} = h_{12} T_{x_1}T^{-1}_{a_1} h_{12}^{-1}
\,\,\, \mbox{and} \,\,\,
T_{x_3}T^{-1}_{a_3} = h_{13} T_{x_1}T^{-1}_{a_1} h_{13}^{-1}
$$
Therefore this group will contain the twist $T_{a_4} = T_{\gamma_k}$, and its
conjugate $T_{\alpha_{k+1}}$. If $I$ also send some $\gamma_i$ into $\beta_j$,
it will imply that the group generated by $\rho_i$ and $I$ contains
$\Mod(\Sigma')$, which allows us to conclude that it is almost the
whole mapping class group.

The lanterns $S$ and $R^2S$ have a common boundary component $a_1=R^2a_2$
and their union is a surface $S_4$ homeomorphic to a  sphere with
$6$ boundary components. There exists an involution $I$ of $S_4$
which takes $S$ to $R^2S$ and acts on the
boundary components as follows
$$
Ia_1 = R^2 a_2 = a_1 \quad Ia_2 = R^2 a_3 \quad
Ia_3 = R^2 a_1 \quad Ia_4 = R^2 a_4
$$
and therefore acts on the middle curves as
$$
Ix_1 = R^2x_2 \quad
Ix_2 = R^2x_3 \quad
Ix_3 = R^2x_1.
$$
Therefore, for any extension $\bar I$ of $I$ to an involution of the surface
$\Sigma$ we have:
$$
T_{x_1}T^{-1}_{a_1} = \rho_1 \rho_3
\quad
T_{x_1}T^{-1}_{a_1} = R^{-2}\bar I\rho_1 \rho_3 \bar I R^2
\quad
T_{x_1}T^{-1}_{a_1} = \bar I R^2\rho_1 \rho_3R^{-2}\bar I.
$$
This shows that the group generated by $\rho_i$ and $\tilde I$ contains
the twist $T_{a_4} = T_{\gamma_k}$ and by conjugation all twists around
$\alpha_i$-es and $\gamma_i$-es. In order to conclude that this group is
the whole mapping class group we need to make sure that $\tilde I$ sends
some $\gamma_i$ to some $\beta_i$ and acts permutes the punctures in
a nice way.

\begin{figure}
\center{
\setlength{\unitlength}{0.00765in}
\begin{picture}(0,0)(0,0)
\put(70,250){${}_{\alpha_{k+1}}$}

\put( 70,125){${}_{\alpha_{k-1}}$}
\put( 95,185){${}_{\gamma_{k-1}}$}
\put( 45,185){${}_{\gamma_{k}}$}

\put( 95,285){${}_{\alpha_{k+3}}$}
\put( 70,350){${}_{\gamma_{k+1}}$}
\put( 45,285){${}_{\gamma_{k+2}}$}

\put(190,305){${}_{\alpha_{k+4}}$}
\put(190,260){${}_{\gamma_{k+3}}$}

\put(190,205){${}_{\beta_{k-2}}$}
\put(190,165){${}_{x}$}

\put(0,245){${}_{I}$}

\put(457,383){${}_{\alpha_{k+1}}$}

\put(250,325){${}_{\alpha_{k-1}}$}
\put(265,340){${}_{\gamma_{k-1}}$}
\put(370,390){${}_{\gamma_{k}}$}

\put(457,275){${}_{\alpha_{k+3}}$}
\put(403,330){${}_{\gamma_{k+1}}$}
\put(403,280){${}_{\gamma_{k+2}}$}

\put(457,220){${}_{\alpha_{k+4}}$}
\put(403,220){${}_{\gamma_{k+3}}$}

\put(315,270){${}_{\beta_{k-2}}$}
\put(335,220){${}_{x}$}

\put(365,445){${}_{\rho_1}$}

\end{picture}
\includegraphics{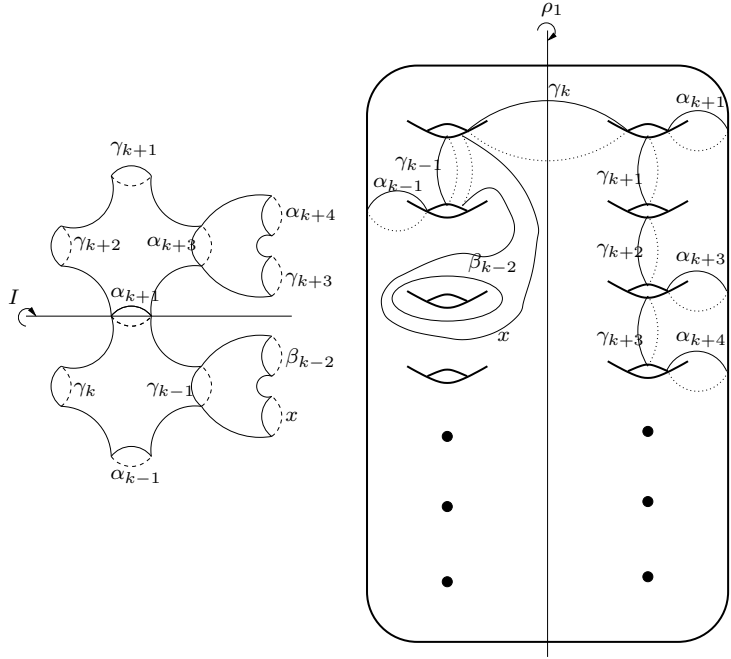}
}
\caption{Extending the lanterns $S$ and $R^2S$ with two pairs of pants.}
\label{3lantern}
\end{figure}

We will extend the surface $S_4$ by attaching two pairs of pants
and obtain a surface $S_5$ homeomorphic to sphere with $8$ boundary components,
such that its compliment $S_6$ in $\Sigma_{g,b}$ is connected.

One of the pair of pants will be bounded by $\alpha_{k+3}=R^2a_1$,
$\gamma_{k+3}$ and $\alpha_{k+4}$, and the second one will have
boundaries $\gamma_{k-1}=a_3 = I(R^2a_1)$, $\beta_{k-2}$ and
some curve $x$ (see Figure \ref{3lantern}). The involution
$I$ can be extended to
$S_5$ with the property  $I\beta_{k-2} = \gamma_{k+3}$.

The compliment $S_6$ is a surface of genus $g-7\geq 2$ with $8$
boundary components and $b$ punctures. Using the arguments from
the previous sections it can be seen that there exists an
involution $\tilde I$ on $S_6$ which swaps the boundary
components in $4$ pairs and acts on the
punctures as the involution $r_3$\footnote{Here we use that we
can find an involution $\tilde I$ of $S_6$ with $s=2(g-7)+2$ fixed
points. The number $s\geq 2$ is greater than the number of fixed points of
$r_3$, which is either $2$ or $3$.}. Let $J$ be the involution
obtained by gluing together $I$ and $\tilde I$.

\begin{theorem}[4 involutions]
\label{4inv}
If $g\geq 7$ for $b$ even and $g>7$ for $b$ odd,
the group generated by $\rho_1$, $\rho_2$, $\rho_3$ and
$J$
is the whole mapping class group $\Mod_{g,b}$.
\end{theorem}
\begin{proof}
Let $G_5$ be the group generated by these $4$ involutions. As explained above,
this group contains the twist $T_{a_k}=T_{\gamma_k}$, and
therefore the twists around all curves $\alpha_i$, $\beta_i$ and $\gamma_i$
because all these curves lie in the orbit of $\gamma_k$ under the group
generated by $I$ and $R$. This allows us to apply Corollary~\ref{genofSb}
and conclude that $G_5$ contains $\Mod^0(\Sigma_{g,b})$.
By construction this group also maps onto $\Sym_b$, because its
image contains the involution $r_1$, $r_2$ and $r_3$, therefore
the group $G_5$ is the whole mapping class group.
\end{proof}
\begin{remark}
The symmetric group $\Sym_b$ can be generated by the involutions $r_1$, $r_2$
and an other involution with $1$ fixed point if $b=4k+3$, which shows that
Theorem~\ref{4inv} can be extended to the case $g=7$ and $b=4k+3$. If $b=4k+1$
this is not possible because every involution with $1$ fixed point is even.
However, we believe that if one slightly modify the action of $\rho_2$ on the
disk $D$ and make $\rho_2$ acts on punctures with $3$ fixed points, it
is possible to obtain that $\Mod_{7,b}$ is generated by $4$ involutions for
all values of $b$.
\end{remark}
\end{subsection}

\begin{subsection}{Genus 3 and odd number of punctures}

Finally we will sketch the case of genus $3$ with odd number of
punctures. In this case it is not possible to extend the pair swap
involutions $J_{12}$ and $J_{23}$ to involutions of the surface
$\Sigma_{3,b}$. Therefore, in order to generate the twist
$T_\alpha$ we need $6$ involutions (one of which is $\rho_1$). We
also need $2$ more involutions which move some curve $\alpha$ to
some curve $\beta$ and $\gamma$, respectively. These, together
with the involution $\rho_2$, give us a generating set consisting
of $9$ involutions.

\begin{theorem}[9 involutions]
\label{9inv}
If $g= 3$ and $b$ is odd, the  mapping class group $\Mod_{g,b}$
can be generated by $9$ involutions.
\end{theorem}
\end{subsection}

This finishes the prof of the last part of Theorem~\ref{main}.


\medskip

\begin{remark}
One can try to find a small generating set of the group
$\Mod^0_{g,b}$ consisting of involutions. It is not clear that
this group is generated by involutions -- if the number of
punctures $b$ is more than $2g+2$ then there are no involutions in
the group $\Mod^0_{g,b}$. It is easy to see that if $g\geq 3$ and
$b \leq 2(g-2)$ then the group $\Mod^0_{g,b}$ is generated by
finitely many involutions. It is interesting to see for which
values of $g$ and $b$ is this group generated by finitely many
involutions and what is the smallest number of involutions needed
to generate this group.

The constructions used in the prof of Theorem~\ref{main} can be
generalized to show that if $b \leq 2(g-2)$, the group
$\Mod^0_{g,b}$ can be generated by $Cg$ involutions for some
constant $C$, and that there is a function $f(\lambda)$ for
$\lambda < 2$ such that if $b \leq \lambda g$ then $\Mod^0_{g,b}$
the group $\Mod^0_{g,b}$ can be generated by $f(\lambda)$
involutions. It is interesting whether there exists a constant $C$
such that the mapping class group $\Mod^0_{g,b}$ can be generated
by $C$ involutions for all values of $g$ and $b$ such that this
group is generated by involutions.
\end{remark}

\end{section}

\bibliographystyle{abbrv}
\bibliography{topology}

\noindent
Martin Kassabov:\\
Department of Mathematics \\
University of Alberta\\
632 Central Academic Building\\
Edmonton, Alberta, T6G 2G1\\
Canada\\
E-mail: kassabov@aya.yale.edu

\end{document}